\documentclass[a4paper, fullpage, reqno, 11pt]{amsart}
\usepackage{amssymb}
\usepackage{ifthen} 

\provideboolean{shownotes} 
\setboolean{shownotes}{true} 
\newcommand{\margnote}[1]{
\ifthenelse{\boolean{shownotes}}%
{\marginpar{\raggedright\tiny\texttt{#1}}}%
{}%
}

\newcommand{\hole}[1]{
\ifthenelse{\boolean{shownotes}}%
{\begin{center} \fbox{ \rule {.25cm}{0cm}
\rule[-.1cm]{0cm}{.4cm} \parbox{.85\textwidth}{\begin{center}
\texttt{#1}\end{center}} \rule {.25cm}{0cm}}\end{center}}
{}
}
\newtheorem{theorem}{Theorem}[section]
\newtheorem{proposition}[theorem]{Proposition}
\newtheorem{lemma}[theorem]{Lemma}
\newtheorem{corollary}[theorem]{Corollary}

\theoremstyle{remark}
\newtheorem{remark}[theorem]{Remark}
\newtheorem{definition}[theorem]{Definition}

\newcommand{\e}{\varepsilon}		       
\newcommand{\R}{\mathbb{R}}

\newcommand{\ue}{u^{\varepsilon}}
\newcommand{\ut}{\tilde{u}}
\newcommand{\pe}{p^{\varepsilon}}
\newcommand{\pt}{\tilde{p}}

\numberwithin{equation}{section}

\begin{document}
\title [Artificial compressibility approximations for Navier Stokes]
{A dispersive approach to the artificial compressibility approximations of the Navier Stokes equations in 3-D}

\author{Donatella Donatelli}
\address{Donatella Donatelli --- 
Dipartimento di Matematica Pura ed Applicata \\
Universit\`a di L'Aquila \\
Via Vetoio, \\
		     67010  Coppito (AQ), Italy}
\email{donatell@univaq.it}

\author{Pierangelo Marcati}
\address{Pierangelo Marcati --- 
Dipartimento di Matematica Pura ed Applicata \\
Universit\`a di L'Aquila \\
Via Vetoio, \\
		     67010  Coppito (AQ), Italy}
\email{marcati@univaq.it}

\begin{abstract}
 In this paper we study how  to approximate the Leray weak solutions of the incompressible Navier Stokes equation. In particular we describe an hyperbolic version of the so called artificial compressibility method investigated by J.L.Lions and Temam. By exploiting the wave equation structure of the pressure  of the approximating system we achieve the convergence of the approximating sequences by means of dispersive  estimate of Strichartz type. We prove that the projection of the approximating  velocity fields on the divergence free vectors is relatively compact and converges to a Leray weak solution of the incompressible Navier Stokes equation. \end{abstract}
\subjclass{35L65}
\keywords{compressible Navier Stokes equation; hyperbolic equations; wave equations.}
\date{}
\maketitle
\section{Introduction}
This paper is concerned with the convergence of the artificial compressibility approximation to the Leray weak solutions (``turbulent in the Leray terminology'') of the $3-D$ Navier Stokes equation on the whole space. This approximation was introduced by Chorin \cite{Ch68, Ch69}, Temam \cite{Tem69a, Tem69b} and Oskolkov \cite{Osk}, in order to deal with the difficulty induced by the incompressibility constraints in the numerical approximations to the Navier Stokes equation. The paper of Temam \cite{Tem69a, Tem69b} and his book \cite{Tem01} discuss the convergence of these approximations on bounded domains by using the classical Sobolev compactness embedding and they recover compactness in time by the classical Lions \cite{L-JL59} method of fractional derivatives. This paper will take a different point of view, namely we wish to exploit the underlying wave equation structure and the presence of dispersive type estimates. In particular we will consider the following system
\begin{equation}
\begin{cases}
\displaystyle{\partial_{t}\ue+\nabla \pe=\mu\Delta \ue-\left(\ue\cdot\nabla\right)\ue -\frac{1}{2}(div \ue)\ue}+f^{\e}\\
\e\partial_{t}\pe+ div \ue=0,
\end{cases}
\label{3}
\end{equation}
where $(x,t)\in \R^{3}\times [0,T]$, $\ue=\ue(x,t)\in\R^{3}$  and $\pe=\pe(x,t)\in \R$, $f^{\e}=f^{\e}(x,t)\in \R^{3}$. \\
The system will be discussed as a semilinear wave type equation for the pressure function and the dispersive estimates will be carried out by using the $L^{p}$-type estimates due to Strichartz \cite{GV95, KT98, S77}. The particular type of  Strichartz estimates that we are going to use here can be found in the book of Sogge \cite{So95} or deduced by the so called bilinear estimates of Klainerman and Machedon \cite{KM93} and  Foschi Klainerman \cite{FK00}. Our analysis can also be related to the convergence of the incompressible limit problem via a formal expansion (see for instance Temam \cite{Tem01}, Chapter 3). In particular a similar wave equation structure has been exploited in various way by the paper of P.L.Lions and Masmoudi \cite{L-P.L.M98}, Desjardin, Grenier, Lions, Masmoudi \cite{DGLM}, Desjardin Grenier \cite{DG99}.\\
In this paper we analyze the convergence problem in the case of the whole space but our method can be extended to exteriors domains which will be done in a forthcoming paper.\\
The interest into the artificial compressibility methods started with the previously mentioned results of Chorin and Temam and was later on investigated by Ghidaglia and Temam \cite{GhT88}. Later developments of numerical investigations in the directions of projections methods have been carried out by \cite{GMS06}, \cite{E03}, \cite{Pro97}, \cite{Ran92}, \cite{NP}, \cite{KS}, \cite{Sma}, \cite{YKS}.\\
This paper is organized as follows. In Section 2 we recall the mathematical tools needed in the paper and recall same basic definitions. In Section 3 we set up our problem, we explain our approximating system and we state our main result. The Section 4 is devoted to recover the a priori estimates needed to  get the strong convergence of the approximating sequences and to prove the main theorem. Finally in Section 5 we give the proof of the main result.

\section{Preliminaries}

For convenience of the reader we establish some notations and recall some basic facts that will be useful in the sequel.\\
We will denote by  $\mathcal{D}(\R ^d \times \R_+)$
 the space of test function
$C^{\infty}_{0}(\R^d \times \R_+)$, by $\mathcal{D}'(\R^d \times
\R_+)$ the space of Schwartz distributions and $\langle \cdot, \cdot \rangle$
the duality bracket between $\mathcal{D}'$ and $\mathcal{D}$ and by $\mathcal{M}_{t}X'$ the space $C_{c}^{0}([0,T];X)'$. Moreover
$W^{k,p}(\R^{d})=(I-\Delta)^{-\frac{k}{2}}L^{p}(\R^{d})$ and $H^{k}(\R^{d})=W^{k,2}(\R^{d})$ denote the nonhomogeneous Sobolev spaces for any $1\leq p\leq \infty$ and $k\in \R$. $\dot W^{k,p}(\R^{d})=(I-\Delta)^{-\frac{k}{2}}L^{p}(\R^{d})$ and $\dot H^{k}(\R^{d})=W^{k,2}(\R^{d})$  denote the homogeneous Sobolev spaces. The notations
 $L^{p}_{t}L^{q}_{x}$ and $L^{p}_{t}W^{k,q}_{x}$ will abbreviate respectively  the spaces $L^{p}([0,T];L^{q}(\R^{d}))$, and $L^{p}([0,T];W^{k,q}(\R^{d}))$.\\
We shall denote by $Q$ and $P$ respectively  the Leray's projectors $Q$ on the space of gradients vector fields and $P$ on the space of divergence - free vector fields. Namely
\begin{equation}
Q=\nabla \Delta^{-1}div\qquad P=I-Q.
\label{1}
\end{equation} 
Let us remark that   $Q$ and $P$  can be expressed in terms of Riesz multipliers, therefore they are  bounded linear operators on every $W^{k,p}$ $(1\leq p\leq\infty)$ space (see \cite{Ste93}).   \\ \\
Let us recall that if  $w$ is a (weak) solution of the following wave equation in the space $[0,T]\times \R^{d}$
\begin{equation*}
\begin{cases}
\left(-\frac{\partial ^{2}}{\partial t}+\Delta\right)w(t,x)=F(t,x)\\
w(0,\cdot)=f,\quad \partial_{t}w(0,\cdot)=g,
\end{cases}
\end{equation*}
for some data $f,g, F$ and time $0<T<\infty$, 
then $w$ satifies the following Strichartz estimates, (see \cite{GV95}, \cite{KT98})
\begin{equation}
\|w\|_{L^{q}_{t}L^{r}_{x}}+\|\partial_{t}w\|_{L^{q}_{t}W^{-1,r}_{x}}\lesssim \|f\|_{\dot H^{\gamma}_{x}}+\|g\|_{\dot H^{\gamma -1}_{x}}+\|F\|_{L^{\tilde{q}'}_{t}L^{\tilde{r}'}_{x}},
\label{s2}
\end{equation}
where $(q,r)$, $(\tilde{q},\tilde{r})$ are \emph{wave admissible} pairs, namely they satisfy 
\begin{equation*}
\frac{2}{q}\leq (d-1)\left(\frac{1}{2}-\frac{1}{r}\right) \qquad 
\frac{2}{\tilde{q}}\leq (d-1)\left(\frac{1}{2}-\frac{1}{\tilde{r}}\right)
\end{equation*}
and moreover the following   conditions holds
\begin{equation*}
\frac{1}{q}+\frac{d}{r}=\frac{d}{2}-\gamma=\frac{1}{\tilde{q}'}+\frac{d}{\tilde{r}'}-2.
\end{equation*}
Later on we shall use \eqref{s2} in the case of $d=3$, $(\tilde{q}', {\tilde{r}'})=(1, 3/2)$,  then $\gamma=1/2$ and $(q,r)=(4,4)$, namely the following estimate
\begin{equation}
\|w\|_{L^{4}_{t,x}}+\|\partial_{t}w\|_{L^{4}_{t}W^{-1,4}_{x}}\lesssim \|f\|_{\dot H^{1/2}_{x}}+\|g\|_{\dot H^{ -1/2}_{x}}+\|F\|_{L^{1}_{t}L^{3/2}_{x}}.
\label{s3}
\end{equation}
Beside the Strichartz estimate \eqref{s2} or \eqref{s3} in the case of $d=3$ (see \cite{So95}), a more refined estimate,related to an earlier linear Strichartz \cite{S77} estimate, can also be deduced by the bilinear estimates of Klainerman and Machedon \cite{KM93}, Foschi  and Klainerman \cite{FK00}, namely
\begin{equation}
\|w\|_{L^{4}_{t,x}}+\|\partial_{t}w\|_{L^{4}_{t}W^{-1,4}_{x}}\lesssim \|f\|_{\dot H^{1/2}_{x}}+\|g\|_{\dot H^{1/2}_{x}}+\|F\|_{L^{1}_{t}L^{2}_{x}}.
\label{s1}
\end{equation}

\section{Approximating system and main result}
Let us consider the incompressible Navier Stokes equation
\begin{equation}
    \begin{cases}
	\partial_{t}u+div(u\otimes u)-\mu\Delta u=\nabla p+f\\
	div u=0\\
	u(x,0)=u_{0},
	\end{cases}
    \label{2}
\end{equation}
where $(x,t)\in \R^{3}\times [0,T]$, $u\in\R^{3}$ denotes the velocity vector field , $p\in \R$ the pressure of 
the fluid , $f\in \R^{3}$ is a given external force, $\mu$ is the kinematic viscosity.
Let us recall (see P.L.Lions \cite{LPL96} and Temam \cite{Tem01}) the notion of Leray weak solution.
\begin{definition}
We say that $u$ is a Leray weak solution of the Navier Stokes equation if it 
satisfies  \eqref{2} in the sense of distributions, namely
\begin{align*}
&\int_{0}^{T}\!\!\int_{\R^{d}}\left(\nabla u\cdot\nabla\varphi -u_{i}u_{j}\partial_{i}\varphi_{j}-u\cdot
\frac{\partial \varphi}{\partial t}\right) dxdt\\&
=\int_{0}^{T}\langle f ,\varphi \rangle_{H^{-1}\times H^{1}_{0}} dxdt+\int_{\R^{d}}u_{0}\cdot \varphi dx,
\end{align*}
for all $\varphi\in C^{\infty}_{0}(\R^{d}\times[0,T])$, $div \varphi =0$ and
\begin{equation*}
div u=0 \qquad \text{in $\mathcal{D}'(\R^{d}\times[0,T])$}
\end{equation*}
and the following energy inequality holds
\begin{align*}
\frac{1}{2}\int_{\R^{d}}&|u(x,t)|^{2}dx+\mu\int_{0}^{t}\!\!\int_{\R^{d}}|\nabla u(x,t)|^{2}dxds\\ \leq
&\frac{1}{2}\int_{\R^{d}}|u_{0}|^{2}dx +\int_{0}^{t}\langle f ,u \rangle_{H^{-1}\times H^{1}_{0}}ds,\qquad \text{for all $t\geq 0$}.
\end{align*}

\end{definition}
There exists in the mathematical literature several results concerning the existence of Leray weak solutions to the Navier Stokes equations, for example we can refer to books of P.L.Lions \cite{LPL96} and Temam \cite{Tem01}. The case $d=3$ is a major open problem and a considerably more difficult case than the case $d=2$, since the bound on the $L^{2}$ norm (kinetic energy) provides only a control on a supercritical norm and does not provide any information concerning the critical controlling (and scaling 
invariant) norm $L^{3}$. Hence we do not know (opposite to the case $d=2$) whether or not the Leray weak solutions are unique, unless (see Serrin \cite{Se63}) we assume a control on the $L^{3}$
norm. Some important regularity results can be found in \cite{CKN}.\\
In order to approximate the system\eqref{2} we wish to use the system  \eqref{3} where we introduce  a ``linearized'' compressibility constraint given by the equation
\begin{equation*}
\partial_{t}\pe=-\frac{1}{\e}div\ue.
\end{equation*}
In order to avoid the paradox of increasing the kinetic energy along the motion we introduce the correction
\begin{equation*}
-\frac{1}{2}(div\ue)\ue
\end{equation*}
into the momentum balance equation.\\
The limiting behaviour as $\e\downarrow 0$ of the initial data to \eqref{3} deserves a little discussion. Indeed \eqref{3} requires two initial conditions
\begin{equation}
\ue(x,0)=\ue_{0}(x), \qquad \pe(x,0)=\pe_{0}(x),
\end{equation}
while the Navier Stokes equations require only one initial condition on the velocity $u$. Hence our approximation will be consistent if the initial datum on the pressure will be eliminated by an ``initial layer'' phenomenon. Since in the limit we have to deal with Leray solutions it is reasonable to require the finite energy constraint to be satisfied by the approximating sequences $(\ue,\pe)$. So we can deduce a natural behaviour to be imposed on the initial data $(\ue_{0},\pe_{0})$, namely
\begin{align}
\tag {\bf{ID}}
& \ue_{0}=\ue(\cdot, 0)\longrightarrow u_{0}=u(\cdot ,0)\ \text{strongly in}\  L^{2}(\R^{3})
\\ &\sqrt{\e}\pe_{0}=\sqrt\e \pe(\cdot,0)\longrightarrow 0\  \text{strongly in} \  L^{2}(\R^{3})\notag.
\end{align}
Let us remark that  the convergence of $\sqrt{\e}\pe_{0}$ to $0$ is necessary to avoid  the presence of concentrations of energy in the limit and it includes the Temam's assumption that $\{\pe_{0}\}$ is bounded in $L^{2}$.\\
Since it will not affect our approximation process, for semplicity from now on, we will take $\mu=1$ and $f^{\e}=0$.  For convenience, let us now formulate an existence theorem concerning the approximating  problem \eqref{3}.
\begin{theorem}
Let $(\ue_{0},\pe_{0})$ satisfy the conditions (ID) for some $\e>0$. Then the system \eqref{3} has a weak solution $(\ue,\pe)$ with the following properties
\begin{itemize}
  \item [\bf{(i)}] $\ue\in L^{\infty}([0,T];L^{2}(\R^{3}))\cap L^{2}([0,T];\dot H^{1}(\R^{3})) $.
  \item [\bf{(ii)}] $\sqrt{\e}\pe \in L^{\infty}([0,T];L^{2}(\R^{3}))$,
\end{itemize}
for all $T>0$.
\end{theorem}
The proof of this theorem will be omitted since it will be a consequence of all the ``a priori bounds'' that will be obtained in the sequel and it will follow from the use of standard finite dimensional Galerkin type approximations.\\
Let us now state our main result. The convergence of $\{\ue\}$ will be  described by analyzing the convergence of the associated Hodge decomposition.
\begin{theorem}
Let $(\ue,\pe)$ be a sequence of weak solution in $\R^{3}$ of the system \eqref{3}, assume that the initial data satisfy (ID). Then 
\begin{itemize}
  \item [\bf{(i)}] There exists $u\in L^{\infty}([0,T];L^{2}(\R^{3}))\cap L^{2}([0,T];\dot H^{1}(\R^{3}))$ such that 
  \begin{equation*}
\ue\rightharpoonup u \quad \text{weakly in $L^{2}([0,T];\dot H^{1}(\R^{3}))$}.
\end{equation*}
  \item [\bf{(ii)}] The gradient component $Q\ue$ of the vector field $\ue$ satisfies
  \begin{equation*}
Q\ue\longrightarrow 0\quad \text{ strongly in $L^{2}([0,T];L^{p}(\R^{3}))$, for any $p\in [4,6)$}.
\end{equation*}
 \item [\bf{(iii)}] The divergence free component $P\ue$ of the vector field $\ue$ satisfies
   \begin{equation*}
P\ue\longrightarrow Pu=u\quad \text{strongly  in $L^{2}([0,T];L^{2}_{loc}(\R^{3}))$}.
\end{equation*}
 \item [\bf{(iv)}] The sequence $\{\pe\}$ will converge in the sense of distribution (more precisely in $H^{-1}_{t}W^{-2,4}_{x}+\mathcal{M}_{t}W^{-1,4/3}_{x}+L^{2}_{t}H^{-1}_{x}$) to 
 \begin{equation*}
p=\Delta^{-1}div \left((u\cdot\nabla)u\right)=\Delta^{-1}tr((Du)^{2}).
\end{equation*}
\end{itemize}
Moreover $u=Pu$ is a Leray weak solution to the incompressible Navier Stokes equation
\begin{equation*}
P(\partial_{t} u-\Delta u+(u\cdot\nabla)u)=0 \quad \text{in $\mathcal{D}'([0,T]\times \R^{3})$},
\end{equation*}
and the following energy inequality holds
\begin{equation}
\frac{1}{2}\int_{\R^{3}}|u(x,t)|^{2}dx+\int_{0}^{T}\!\!\int_{\R^{3}}|\nabla u(x,t)|^{2}dxdt\leq 
\frac{1}{2}\int_{\R^{3}}|u(x,0)|^{2}dx.
\label{en}
\end{equation}
\label{tM}
\end{theorem}
\begin{remark}
This theorem can be easily extended to the nonhomegeneous equation \eqref{2}, by assuming 
\begin{equation*}
f^{\e}\longrightarrow f\qquad \text{strongly in $L^{2}([0,T];H^{-1}(\R^{3}))$}.
\end{equation*}
\end{remark}
\begin{remark}
Let us denote by $R_{j}$ the Riesz transform. The Hardy space $\mathcal{H}^{1}(\R^{3})$ is a closed subspace of $L^{1}(\R^{3})$ defined by
\begin{equation*}
\mathcal{H}^{1}(\R^{3})=\{f\in L^{1}(\R^{3})\mid R_{j}f\in  L^{1}(\R^{3}), \text{for any}\  j=1,\ldots 3\}.
\end{equation*}
Then one has
\begin{equation}
p\in L^{1}([0,T];L^{3/2}(\R^{3}))\cap L^{1}([0,T];L^{3}(\R^{3})),
\end{equation}
and there exits $c_{1}>0$, such that
\begin{equation*}
\|(tr(Du)^{2})\|_{L^{1}([0,T];\mathcal{H}^{1}(\R^{3}))}\leq c_{1}\|u_{0}\|_{L^{2}(\R^{3})}^{2}.
\end{equation*}

\end{remark}
\section{A priori estimates}
In this section we wish to establish the priori estimates, independent on $\e$, for the solutions of the system \eqref{3} which are necessary to prove the Theorem \ref{tM}.
 We will achieve this goal in two steps. First of all we will recover the  a priori estimates that come from the classical energy estimates related to the system \eqref{3}. Then we get  stronger  estimates by exploiting the structure of the system. 
In fact, as we will see later on , the sequence $\pe$ satifies a  wave type equation. This will allow us to apply to $\pe$ the Strichartz estimates \eqref{s1}, \eqref{s2}, and to get in this way dispersive  bounds on $\pe$.
 \subsection{Energy estimates}
The next results concerns the  energy type estimate for the system \eqref{3}.
\begin{theorem}
Let us consider the solution $(\ue, \pe)$ of the Cauchy problem for the system \eqref{3}. Assume that the hypotheses (ID) hold, then one has
\begin{equation}
\label{9}
E(t)+\int_{0}^{t}\!\!\int_{\R^{3}}|\nabla \ue(x,s)|^{2}dxds=E(0),
\end{equation}
where we set
\begin{equation}
\label{8}
E(t)=\int_{\R^{3}}\left( \frac{1}{2}|\ue(x,t)|^{2}+ \frac{\e}{2} |\pe(x,t)|^{2}\right)dx.
\end{equation}
\label{t1}
\end{theorem}
\begin{proof}
We multiply, as usual,  the first equation of the system \eqref{3} by $\ue$ and the second by $\pe$, then we sum up and integrate by parts in space and time, hence  we get \eqref{9}.
\end{proof}
\begin{corollary}
\label{c1}
Let us consider the solution $(\ue, \pe)$ of the Cauchy problem for the system \eqref{3}. Let us assume that the hypotheses (ID) hold, then  it follows
\begin{align} 
& \sqrt{\e}\pe &\quad & \text{is bounded in $L^{\infty}([0,T];L^{2}(\R^{3}))$,} \label{11}\\
& \e\pe_{t} &\quad & \text{is relatively compact in $H^{-1}([0,T]\times \R^{3}),$}  \label{12}\\
& \nabla\ue &\quad & \text{is bounded in $L^{2}([0,T]\times\R^{3}),$}  \label{13}\\
& \ue &\quad & \text{is bounded in $L^{\infty}([0,T];L^{2}(\R^{3}))\cap L^{2}([0,T];L^{6}(\R^{3})),$}  \label{14}\\
(&\ue \!\cdot\!\nabla)\ue &\quad & \text{is bounded in $L^{2}([0,T];L^{1}(\R^{3}))\cap L^{1}([0,T];L^{3/2}(\R^{3})),$}  \label{17}\\
(& div\ue)\ue &\quad & \text{is bounded in $L^{2}([0,T];L^{1}(\R^{3}))\cap L^{1}([0,T];L^{3/2}(\R^{3})).$}\label{17a}
\end{align}
\end{corollary}
\begin{proof}
\eqref{11},  \eqref{12}, \eqref{13} follow from \eqref{9}, while \eqref{14} follows from \eqref{9} and Sobolev embeddings theorems. Finally \eqref{17} and \eqref{17a} come from \eqref{13}, \eqref{14}. 
\end{proof}
\subsection{Pressure wave equation}
In this section by using the Strichartz estimates \eqref{s3}, \eqref{s1} we get a priori estimates on $\pe$. We will use a wave equation structure for $\pe$. First of all let us rescale the time variable, the velocity and the pressure in the following way 
\begin{equation}
\label{18}
\tau=\frac{t}{\sqrt{\e}}, \quad \ut(x,\tau)=\ue(x,\sqrt{\e}\tau), \quad \pt(x,\tau)=\pe(x,\sqrt{\e}\tau).
\end{equation}
As a consequence of this scaling the system \eqref{3} becomes
\begin{equation}
\begin{cases}
\displaystyle{\partial_{\tau}\ut+\sqrt{\e}\nabla \pt=\sqrt{\e}\Delta \ut-\sqrt{\e}\left(\ut\cdot\nabla\right)\ut -\frac{\sqrt{\e}}{2}(div \ut)\ut}\\
\sqrt{\e}\partial_{\tau}\pt+ div \ut=0
\end{cases}
\label{19}
\end{equation}
then, by differentiating with respect to time the equation $\eqref{19}_{2}$ and by using $\eqref{19}_{1}$, we get that $\pt$ satisfies the following wave equation
\begin{equation}
\label{20 }
\partial_{\tau\tau}\pt-\Delta\pt +\Delta div\ut-div\left(\left(\ut\cdot\nabla\right)\ut +\frac{1}{2}(div \ut)\ut \right)=0.
\end{equation}
Now we consider $\pt=\pt_{1}+\pt_{2}$ where $\pt_{1}$ and $\pt_{2}$  solve the following wave equations:
\begin{equation}
\label{21}
\begin{cases}
     \partial_{\tau\tau}\pt_{1}-\Delta\pt_{1} =-\Delta div\ut=F_{1} \\
     \pt_{1}(x,0)=\partial_{\tau} \pt_{1}(x,0)=0,
   \end{cases}
\end{equation}    
\begin{equation}
\label{22} 
\begin{cases}   
\displaystyle{ \partial_{\tau\tau}\pt_{2}-\Delta\pt_{2} =div\left(\left(\ut\cdot\nabla\right)\ut +\frac{1}{2}(div \ut)\ut \right)=F_{2}}\\
\pt_{2}(x,0)=\pt(x,0)\quad  \partial_{\tau}\pt_{2}(x,0)=\partial_{\tau}\pt(x,0).
 \end{cases}
 \end{equation}
Therefore we are able to prove the following theorem.
\begin{theorem}
Let us consider the solution $(\ue, \pe)$ of the Cauchy problem for the system \eqref{3}. Assume that the hypotheses (ID) hold. Then we set the following estimate 
\begin{align}
\hspace{-1mm}\e^{3/8}\|\pe\|_{L^{4}_{t} W^{-2,4}_{t}}+\e^{7/8}\|\partial_{t}\pe\|_{L^{4}_{t} W^{-3,4}_{t}}&\lesssim \sqrt{\e}\|\pe_{0}\|_{L^{2}_{x}}+\|div\ue_{0}\|_{H^{-1}_{x}}+\sqrt{T}\|div \ue\|_{L^{2}_{t}L^{2}_{x}}\notag\\&+
\|\left(\ue\cdot\nabla\right)\ue +\frac{1}{2}(div \ue)\ue\|_{L^{1}_{t}L^{3/2}_{x}}.
\label{23}
\end{align}
\label{t2}
\end{theorem}
\begin{proof}
Since $\pt_{1}$ and $\pt_{2}$ are solutions of the wave equations \eqref{21}, \eqref{22}, we can apply the Strichartz estimates \eqref{s3} and \eqref{s1}, with $(x,\tau)\in \R^{3}\times\left (0,T/\sqrt \e\right)$.
Since  $\Delta ^{-1}\pt_{1}$ satisfies the equation
\begin{equation}
\label{ 24}
 \partial_{\tau\tau}(\Delta^{-1}\pt_{1})-\Delta(\Delta^{-1}\pt_{1}) =\Delta ^{-1}F_{1},
\end{equation}
then by using the Strichartz estimates \eqref{s1} we get
\begin{equation}
\label{25}
\|\Delta^{-1}\pt_{1}\|_{L^{4}_{\tau,x}}+\|\partial_{\tau}\Delta^{-1}\pt_{1}\|_{L^{4}_{\tau} W^{-1,4}_{x}}\lesssim 
\|\Delta^{-1}F_{1}\|_{L^{1}_{\tau}L^{2}_{x},}
\end{equation}
namely
\begin{equation}
\|\pt_{1}\|_{L^{4}_{\tau} W^{-2,4}_{x}}+\|\partial_{\tau}\pt_{1}\|_{L^{4}_{\tau} W^{-3,4}_{x}}\lesssim \frac{\sqrt{T}}{\e^{1/4}}\|div \ut\|_{L^{2}_{\tau}L^{2}_{x}}.
\label{26}
\end{equation}
In the same way we have that $\Delta ^{-1/2}\pt_{2}$ satisfies the equation
\begin{equation}
\label{ 27}
 \partial_{\tau\tau}(\Delta^{-1/2}\pt_{2})-\Delta(\Delta^{-1/2}\pt_{1}) =\Delta ^{-1/2}F_{2},
\end{equation}
therefore by using  the estimate \eqref{s3} we get
\begin{align}
\label{28}
\|\Delta^{-1/2}\pt_{2}\|_{L^{4}_{\tau,x}}+\|\partial_{\tau}\Delta^{-1/2}\pt_{2}\|_{L^{4}_{\tau} W^{-1,4}}&\lesssim \|\Delta^{-1/2}\pt(x,0)\|_{ H^{1/2}_{x}}\notag\\&+
\|\Delta^{-1/2}\partial_{\tau}\pt(x,0)\|_{ H^{-1/2}_{x}}\notag\\&+\|\Delta^{-1/2}F_{2}\|_{L^{1}_{\tau}L^{3/2}_{x},}
\end{align}
namely
\begin{align}
\label{29}
\|\pt_{2}\|_{L^{4}_{\tau} W^{-1,4}_{x}}+\|\partial_{\tau}\pt_{2}\|_{L^{4}_{\tau} W^{-2,4}_{x}}&\lesssim \|\pt(x,0)\|_{ H^{-1/2}_{x}}+
\|\partial_{\tau}\pt(x,0)\|_{ H^{-3/2}_{x}}\notag
\\&+\|\left(\ut\cdot\nabla\right)\ut +\frac{1}{2}(div \ut)\ut \|_{L^{1}_{\tau}L^{3/2}_{x},}
\end{align}
Now by using \eqref{26}, \eqref{29} it follows that $\pt$ verifies
\begin{align}
\|\pt\|_{L^{4}_{\tau} W^{-2,4}_{x}}+\|\partial_{\tau}\pt\|_{L^{4}_{\tau} W^{-3,4}_{x}}&\leq \|\pt_{1}\|_{L^{4}_{\tau} W^{-2,4}_{x}}+\|\pt_{2}\|_{L^{4}_{\tau} W^{-1,4}_{x}}\\&+
\|\partial_{\tau}\pt_{1}\|_{L^{4}_{\tau} W^{-3,4}_{x}}+\|\partial_{\tau}\pt_{2}\|_{L^{4}_{\tau} W^{-2,4}_{x}}
\notag
\\&\lesssim \|\pt(x,0)\|_{ H^{-1/2}_{x}}+
\|\partial_{\tau}\pt(x,0)\|_{ H^{-3/2}_{x}}\notag\\&
+\frac{\sqrt{T}}{\e^{1/4}}\|div \ut\|_{L^{2}_{\tau}L^{2}_{x},}+
\|\left(\ut\cdot\nabla\right)\ut +\frac{1}{2}(div \ut)\ut \|_{L^{1}_{\tau}L^{3/2}_{x}}.\notag
\label{30}
\end{align}
Finally, since
\begin{equation*}
\|\pt\|_{L^{r}((0,T/\sqrt \e );L^{q}(\R^{3}))}=\e^{-1/2r}\|\pe\|_{L^{r}([0,T];L^{q}(\R^{3}))}
\end{equation*}
 we end up with \eqref{23}.
\end{proof}
\section{Strong convergence}
In this section we conlcude the proof of the  Theorem \ref{tM}. In particular we will show that the gradient part of the velocity $Q\ue$ converges strongly to $0$, while the incompressible component of the velocity field $P\ue$ converges strongly to $Pu=u$, where $u$ is the limit profile as $\e\downarrow 0$ of $\ue$.
\subsection{Strong convergence of $Q\ue$ and $P\ue$}
We start this section with some easy consequences of the a priori estimates established in the previous section.
\begin{corollary}
\label{c2}
Let us consider the solution $(\ue, \pe)$ of the Cauchy problem for the system \eqref{3}. Assume that the hypotheses (ID) hold. Then, as $\e\downarrow 0$, one has
\begin{align}
&\e\pe\longrightarrow 0 &\quad& \text{strongly in $L^{\infty}([0,T];L^{2}(\R^{3}))\cap L^{4}([0,T];W^{-2,4}(\R^{3}))$,}\label{31}\\
&div \ue \longrightarrow 0 &\quad& \text{strongly in $ W^{-1,\infty}([0,T];L^{2}(\R^{3}))\cap L^{4}([0,T];W^{-3,4}(\R^{3}))$}.\label{32}
\end{align}
\end{corollary}
\begin{proof}
\eqref{31}, \eqref{32} follow  from the estimates \eqref{11}, \eqref{23} and the second equation of the system \eqref{3}. 
\end{proof}
Now, we wish to show that the gradient part of the velocity field $Q\ue$ goes strongly to $0$ as $\e\downarrow 0$. As we will see in the next proposition, this will be a consequence of the estimate \eqref{23} and of the following auxiliary result.
\begin{lemma}
Let us consider  a smoothing kernel $\psi\in C^{\infty}_{0}(\R^{d})$, such that $\psi\geq 0$, $\int_{\R^{d}}\psi dx=1$, and define
\begin{equation*}
\psi_{\alpha}(x)=\alpha^{-d}\psi\left(\frac{x}{\alpha}\right).
\end{equation*}
Then  for any $f\in \dot H^{1}(\R^{d})$, one has
\begin{equation}
\label{y1}
\|f-f\ast \psi_{\alpha}\|_{L^{p}(\R^{d})}\leq C_{p}\alpha^{1-\sigma}\|\nabla f\|_{L^{2}(\R^{d})},
\end{equation}
where
\begin{equation*}
p\in [2, \infty)
\quad \text{if $d=2$}, \quad p\in [2, 6] \quad \text{if $d=3$ \ and}\quad \sigma=d\left(\frac{1}{2}-\frac{1}{p}\right).
\end{equation*}
Moreover the following Young type inequality hold
\begin{equation}
\label{y2}
\|f\ast\psi_{\alpha}\|_{L^{p}(\R^{d})}\leq C\alpha^{s-d\left(\frac{1}{q}-\frac{1}{p}\right)}\|f\|_{W^{-s,q}(\R^{d})},
\end{equation}
for any $p,q\in [1, \infty]$, $q\leq p$,  $s\geq 0$, $\alpha\in(0,1)$.
\label{ly}
\end{lemma}
\begin{proposition}
Let us consider the solution $(\ue, \pe)$ of the Cauchy problem for the system \eqref{3}. Assume that the hypotheses (ID) hold. Then  as $\e\downarrow 0$,
\begin{equation}
Q\ue \longrightarrow 0 \quad \text{strongly in $ L^{2}([0,T];L^{p}(\R^{3}))$ for any $p\in [4,6)$ }.
\label{33}
\end{equation}
\label{p2}
\end{proposition}
\begin{proof}
In order to prove the Proposition \ref{p2} we split $Q\ue$ as follows
\begin{equation*}
\|Q\ue\|_{L^{2}_{t}L^{p}_{x}}\leq \|Q\ue-Q\ue\ast \psi_{\alpha}\|_{L^{2}_{t}L^{p}_{x}}+\|Q\ue\ast \psi_{\alpha}\|_{L^{2}_{t}L^{p}_{x}}=J_{1}+J_{2},
\end{equation*}
where $\psi_{\alpha}$ is the smoothing kernel defined in Lemma \ref{ly}.
Now we estimate separately $J_{1}$ and $J_{2}$. For $J_{1}$ by using \eqref{y1} we get
\begin{equation}
\label{50}
J_{1}\leq \alpha^{1-3\left(\frac{1}{2}-\frac{1}{p}\right)}\left(\int_{0}^{T}\|\nabla Q\ue(t)\|_{L^{2}_{x}}^{2} dt\right)\leq \alpha^{1-3\left(\frac{1}{2}-\frac{1}{p}\right)}\|\nabla \ue\|_{L^{2}_{t}L^{2}_{x}}.
\end{equation}
Hence from the identity $Q\ue=-\e^{1/8}\nabla\Delta^{-1}\e^{7/8}\partial_{t}p$ and by the inequality \eqref{y2} we get  $J_{2}$ satisfies the following estimate
\begin{align}
J_{2}&\leq \e^{1/8}\|\nabla\Delta^{-1}\e^{7/8}\partial_{t}p\ast\psi\|_{L^{2}_{t}L^{p}_{x}}
\leq \e^{1/8}\alpha^{-2-3\left(\frac{1}{4}-\frac{1}{p}\right)}\|\e^{7/8}\partial_{t}p\|_{L^{2}_{t}W^{-3,4}_{x}}\notag\\
&\leq \e^{1/8}\alpha^{-2-3\left(\frac{1}{4}-\frac{1}{p}\right)}T^{1/4}\|\e^{7/8}\partial_{t}p\|_{L^{4}_{t}W^{-3,4}_{x}}.
\label{51}
\end{align}
Therefore,  summing up \eqref{50} and \eqref{51} and by using \eqref{13} and \eqref{23}, we conclude  for any $p\in [4,6)$ that
\begin{equation}
\|Q\ue\|_{L^{2}_{t}L^{p}_{x}}\leq C\alpha^{1-3\left(\frac{1}{2}-\frac{1}{p}\right)}+C_{T}\e^{1/8}\alpha^{-2-3\left(\frac{1}{4}-\frac{1}{p}\right)}.
\end{equation}
Finally we choose $\alpha$ in terms of $\e$ in order that the two terms in the right hand side of the previous inequality have the same order, namely
\begin{equation}
\alpha=\e^{1/18}.
\end{equation}
Therefore we obtain
\begin{equation*}
\displaystyle{\|Q\ue\|_{L^{2}_{t}L^{p}_{x}}\leq C_{T}\e^{ \frac{6-p}{36p}}\quad \text{for any $p\in [4,6)$.}}
\end{equation*}
\end{proof}
It remains to prove the strong compactness of the incompressible component of the velocity field.  To achieve this goal we need to recall here, the following theorem (see \cite{Si}).
\begin{theorem}
Let be $\mathcal{F}\subset L^{p}([0,T];B)$,  $1\leq p<\infty$, $B$ a Banach space. $\mathcal{F}$ is relatively compact in  $L^{p}([0,T];B)$ for $1\leq p<\infty$, or in $C([0,T];B)$ for $p=\infty$ if and only if 
\begin{itemize}
\item[{\bf (i)}]
$\displaystyle{\left\{\int_{t_{1}}^{t_{2}}f(t)dt,\ f\in B\right\}}$ is relatively compact in $B$, $0<t_{1}<t_{2}<T$,
\item[{\bf (ii)}]
$\displaystyle{\lim_{h\to 0}\|f(x+h) - f(x)\|_{L^{p}([0, T-h];B)}=0}$ uniformly for any $f \in \mathcal{F}$.
\end{itemize}
\label{S}
\end{theorem}
The compactness can be obtained by looking at some time regularity properties of $P\ue$ and by using the Theorem \ref{S}, but before we need to prove  the following lemma.
\begin{lemma}
\label{l1}
Let us consider the solution $(\ue, \pe)$ of the Cauchy problem for the system \eqref{3}. Assume that the hypotheses (ID) hold. Then for all $h\in(0,1)$, we have
\begin{equation}
\label{34}
\|P\ue(t+h)-P\ue(t)\|_{L^{2}([0,T]\times \R^{3})}\leq C_{T}h^{1/5}.
\end{equation}
\end{lemma}
\begin{proof}
Let us set $z^{\e}=\ue(t+h)-\ue(t)$, we have
\begin{align}
\hspace{-0,25 cm}\|P\ue(t+h)-P\ue(t)\|^{2}_{L^{2}([0,T]\times \R^{3})}&=\int_{0}^{T}\!\!\int_{\R^{3}}dtdx(Pz^{\e})\cdot(Pz^{\e}-Pz^{\e}\ast \psi_{\alpha})\notag\\&+\int_{0}^{T}\!\!\int_{\R^{3}}dtdx(Pz^{\e})\cdot(Pz^{\e}\ast \psi_{\alpha})=I_{1}+I_{2}.
\label{ 35}
\end{align}
By using \eqref{y1} we can estimate $I_{1}$ in the following way
\begin{align}
I_{1}&\leq \|Pz^{\e}\|_{L^{\infty}_{t}L^{2}_{x}}\int_{0}^{T}\|Pz^{\e}(t)-(Pz^{\e}\ast \psi_{\alpha})(t)\|_{L^{2}_{x}}dt\notag\\&\lesssim \alpha T^{1/2}\|\ue\|_{L^{\infty}_{t}L^{2}_{x}}\|\nabla\ue\|_{L^{2}_{t,x}}.
\label{36}
\end{align}
Let us reformulate $Pz^{\e}$ in integral form by using the equation $\eqref{3}_{1}$, hence
\begin{align}
\hspace{-0.3cm}I_{2}\leq\left|\int_{0}^{T}\!\!\!dt\!\!\int_{\R^{3}}\!\!\!dx \!\!\int_{t}^{t+h}\!\!\!ds(\Delta \ue-\left(\ue\cdot\nabla\right)\ue -\frac{1}{2}\ue(div \ue)(s,x)\cdot (Pz^{\e}\ast \psi_{\alpha})(t,x)\right|.
\label{37}
\end{align}
Then integrating by parts and by using \eqref{y2}, with $p=\infty$ and $q=2$, we deduce
\begin{align}
I_{2}&\leq h\|\nabla\ue\|^{2}_{L^{2}_{t,x}}+C\alpha^{-3/2}T^{1/2}\|\ue\|_{L^{\infty}_{t}L^{2}_{x}}\left(\!h\!\int_{t}^{t+h}\!\!\!\|\left(\ue\cdot\nabla\right)\ue -\frac{1}{2}(div \ue)\ue\|^{2}_{L^{1}_{x}}ds\right)^{1/2}\notag\\
&\leq h\|\nabla\ue\|^{2}_{L^{2}_{t,x}}+C\alpha^{-3/2}T^{1/2}h\|\ue\|_{L^{\infty}_{t}L^{2}_{x}}\|\left(\ue\cdot\nabla\right)\ue -\frac{1}{2}(div \ue)\ue\|_{L^{2}_{t}L^{1}_{x}}.
\label{38}
\end{align}
Summing up $I_{1}$, $I_{2}$ and by taking into account \eqref{13}, \eqref{14}, \eqref{17}, \eqref{17a}, we have
\begin{equation}
\|P\ue(t+h)-P\ue(t)\|^{2}_{L^{2}([0,T]\times \R^{3})}\leq C(\alpha T^{1/2}+h\alpha^{-3/2}T^{1/2}+h),
\label{39}
\end{equation}
by choosing $\alpha=h^{2/5}$, we end up with \eqref{34}.
\end{proof}
\begin{corollary}
Let us consider the solution $(\ue, \pe)$ of the Cauchy problem for the system \eqref{3}. Assume that the hypotheses (ID) hold. Then  as $\e\downarrow 0$
\label{c3}
\begin{equation}
P\ue \longrightarrow Pu, \qquad \text{strongly in $L^{2}(0,T;L^{2}_{loc}(\R^{3}))$}.
\label{43}
\end{equation}
\end{corollary}
\begin{proof}
By using the Lemma \ref{l1} and the Theorem \ref{S} and the Proposition \ref{p2} we get \eqref{43}. 
\end{proof}
\subsection{Proof of the Theorem \ref{tM}}
{\bf (i)} It follows from the estimate \eqref{14}.\\
{\bf (ii)} It is a consequence of  the Proposition \ref{p2}.\\
{\bf (iii)} By taking into account the decomposition $\ue=P\ue+Q\ue$, by the Corollary \ref{c3} and the Proposition \ref{p2} we have that
\begin{equation*}
P\ue\longrightarrow u \qquad \text{strongly in $L^{2}([0,T];L^{2}_{loc}(\R^{3}))$.}
\end{equation*}\\
{\bf (iv)}  Let us apply the Leray projector $Q$ to the equation $\eqref{3}_{1}$, then it follows
\begin{equation}
\label{54}
\nabla \pe =\Delta Q\ue- Q\left(div(\ue\otimes\ue) +\frac{3}{2} \ue div Q\ue\right).
\end{equation}
Now by choosing a test function $\varphi \in H^{1}_{t}W^{2,4/3}_{x}\cap C^{0}_{t}W^{1,4}_{x}\cap L^{2}_{t}H^{1}_{x}$ and by taking into account \eqref{13}, \eqref{33}, \eqref{43}, we get, as $\e \downarrow 0$, 
\begin{align}
\langle\ue div Q\ue, Q\varphi\rangle&\leq\|Q\ue\|_{L^{2}_{t}L^{4}_{x}} \|\nabla\ue\|_{L^{2}_{t}L^{2}_{x}}\|Q\varphi\|_{L^{\infty}_{t}L^{4}_{x}} \notag\\&+ \|Q\ue\|_{L^{2}_{t}L^{4}_{x}} \|\ue\|_{L^{\infty}_{t}L^{2}_{x}}\|\nabla Q\varphi\|_{L^{2}_{t}L^{4}_{x}} \rightarrow 0, 
\end{align}
\begin{align}
\langle div(\ue\otimes\ue),Q\varphi\rangle&=\langle div(P\ue\otimes P\ue),Q\varphi\rangle+\langle div(Q\ue\otimes Q\ue),Q\varphi\rangle\notag\\&+\langle div(P\ue\otimes Q\ue),Q\varphi\rangle\notag+\langle div(Q\ue\otimes Q\ue),Q\varphi \rangle\notag\\&\rightarrow \langle div(Pu \otimes Pu),Q\varphi\rangle=\langle Qdiv((Pu\cdot\nabla)Pu) , \varphi\rangle.
\end{align}
So  as $\e \downarrow 0$ we have, 
\begin{equation}
\label{40}
\langle \nabla \pe , \varphi \rangle \longrightarrow\langle \nabla\Delta^{-1}div((u\cdot\nabla)u) , \varphi\rangle. 
\end{equation}
Now we can pass  into the limit inside the system \eqref{3} and we get  $u$ satisfies the following equation in $\mathcal{D}'([0,T]\times \R^{3})$
\begin{equation}
\label{55}
P(\partial_{t} u-\Delta u+(u\cdot\nabla)u)=0.
\end{equation}
Finally we prove the energy inequality. By using the weak lower semicontinuity of the weak limits, the hypotheses (ID) and denoting 
by $\chi$ the weak-limit of $\sqrt{\e}\pe $,  we have
\begin{align}
&\int_{\R^{3}}\frac{1}{2}|\chi|^{2}dx+\int_{\R^{3}}\frac{1}{2}|u(x,t)|^{2}dx+\int_{0}^{T}\!\!\int_{\R^{3}}|\nabla u(x,t)|^{2}dxdt\notag\\&\leq
\liminf_{\e\to 0}\left(\int_{\R^{3}}\frac{1}{2}|\ue(x,t)|^{2}dx+\int_{\R^{3}}\frac{\e}{2}|\pe|^{2}+\int_{0}^{T}\!\!\int_{\R^{3}}|\nabla \ue(x,t)|^{2}dxdt\right)\notag\\&=\liminf_{\e\to 0}\int_{\R^{3}}\frac{1}{2}\left(|\ue_{0}|^{2}-\e|\pe_{0}|^{2}\right)dx=\int_{\R^{3}}\frac{1}{2}|u_{0}|^{2}dx.
\end{align}


\begin{thebibliography}{0}

\bibitem{CKN}
L.~Caffarelli, R.~Kohn, and L.~Nirenberg, Partial regularity of suitable
  weak solutions of the {N}avier-{S}tokes equations, \emph{Comm. Pure Appl. Math.}
  \textbf{35} (1982), no.~6, 771--831.
  
\bibitem{Ch68}
A.~J.~Chorin, Numerical solution of the {N}avier-{S}tokes
  equations, \emph{Math. Comp.} \textbf{22} (1968), 745--762.

\bibitem{Ch69}
A.~J.~Chorin, On the convergence of discrete approximations to the
  {N}avier-{S}tokes equations, \emph{Math. Comp}. \textbf{23} (1969), 341--353.  
  

\bibitem{DGLM}
B.~Desjardins, E.~Grenier, P.-L. Lions, and N.~Masmoudi, Incompressible
  limit for solutions of the isentropic {N}avier-{S}tokes equations with
  {D}irichlet boundary conditions,\emph{ J. Math. Pures Appl}. (9) \textbf{78}
  (1999), no.~5, 461--471.

\bibitem{DG99}
B.~Desjardins and E.~Grenier, Low {M}ach number limit of viscous
  compressible flows in the whole space, \emph{R. Soc. Lond. Proc. Ser. A Math.
  Phys. Eng. Sci.} \textbf{455} (1999), no.~1986, 2271--2279.

\bibitem{DM04}
D.~Donatelli and P.~Marcati, Convergence of singular limits for multi-{D}
  semilinear hyperbolic systems to parabolic systems, \emph{Trans. Amer. Math. Soc}.
  \textbf{356} (2004), no.~5, 2093--2121 (electronic).

\bibitem{E03}
 E. ~W.  and J.-G. Liu, Gauge method for viscous incompressible
  flows, \emph{Commun. Math. Sci}. \textbf{1} (2003), no.~2, 317--332.
  
\bibitem{FK00}
D.~Foschi and S.~Klainerman, Bilinear space-time estimates for
  homogeneous wave equations, \emph{Ann. Sci. \'Ecole Norm. Sup. (4)} \textbf{33}
  (2000), no.~2, 211--274.
  
\bibitem{GhT88}
J.-M. Ghidaglia and R.~Temam, Long time behavior for partly dissipative
  equations: the slightly compressible {$2$}{D}-{N}avier-{S}tokes equations,
 \emph{ Asymptotic Anal}. \textbf{1} (1988), no.~1, 23--49.

\bibitem{GV95}
J.~Ginibre and G.~Velo, Generalized {S}trichartz inequalities for the
  wave equation, \emph{J. Funct. Anal.} \textbf{133} (1995), no.~1, 50--68.

\bibitem{GMS06}
J.~L. Guermond, P.~Minev, and J.~Shen, An overview of projection methods
  for incompressible flows, \emph{Comp. Meth. Appl. Mech. and Eng}. (2006), to appear.
  
\bibitem{KS}
B.~G. Kuznecov and {\v{S}}.~Smagulov, Approximation of the
  {N}avier-{S}tokes equations, \emph{\v Cisl. Metody Meh. Splo\v sn. Sredy}
  \textbf{6} (1975), no.~2, 70--79.
  
  
\bibitem{KT98}
M.~Keel and T.~Tao, Endpoint {S}trichartz estimates, emph{Amer. J.
  Math.} \textbf{120} (1998), no.~5, 955--980.

\bibitem{KM93}
S.~Klainerman and M.~Machedon, Space-time estimates for null forms and
  the local existence theorem, \emph{Comm. Pure Appl. Math}. \textbf{46} (1993),
  no.~9, 1221--1268.
  
\bibitem{L-JL69}
J.-L. Lions, \emph{Quelques m\'ethodes de r\'esolution des probl\`emes aux
  limites non lin\'eaires}, Dunod, 1969.

\bibitem{L-JL78}
J.-L. Lions, On some problems connected with {N}avier-{S}tokes equations,
  Nonlinear evolution equations \emph{(Proc. Sympos., Univ. Wisconsin, Madison, Wis.,
  1977),} Publ. Math. Res. Center Univ. Wisconsin, vol.~40, Academic Press, New
  York, 1978, pp.~59--84.

\bibitem{L-JL59}
J.-L. Lions, Sur l'existence de solutions des \'equations de
  {N}avier-{S}tokes, \emph{C. R. Acad. Sci. Paris} \textbf{248} (1959), 2847--2849.
 
 \bibitem{L-P.L.M98}
P.-L. Lions and N.~Masmoudi,Incompressible limit for a viscous
  compressible fluid,  \emph{J. Math. Pures Appl}. (9) \textbf{77} (1998), no.~6,
  585--627. 
\bibitem{LPL96}
P.L. Lions, \emph{Mathematical topics in fluid dynamics, incompressible
  models}, Claredon Press, Oxford Science Publications, 1996.  

\bibitem{NP}
R~H. Nochetto and J.-H.~ Pyo, Error estimates for semi-discrete
  gauge methods for the {N}avier-{S}tokes equations, \emph{Math. Comp.} \textbf{74}
  (2005), no.~250, 521--542 (electronic).

\bibitem{Osk}
A.~P. Oskolkov, A certain quasilinear parabolic system with small
  parameter that approximates a system of {N}avier-{S}tokes equations, \emph{Zap.
  Nau\v cn. Sem. Leningrad. Otdel. Mat. Inst. Steklov. (LOMI)} \textbf{21}
  (1971), 79--103.

\bibitem{Pro97}
A.~Prohl, \emph{Projection and quasi-compressibility methods for solving
  the incompressible {N}avier-{S}tokes equations}, Advances in Numerical
  Mathematics, B. G. Teubner, Stuttgart, 1997.

\bibitem{Ran92}
R.~ Rannacher,On {C}horin's projection method for the incompressible
  {N}avier-{S}tokes equations,  The Navier-Stokes equations II---theory and
  numerical methods (Oberwolfach, 1991), \emph{Lecture Notes in Math.}, vol. 1530,
  Springer, Berlin, 1992, pp.~167--183.
  
\bibitem{Se63}
J.~Serrin,The initial value problem for the {N}avier-{S}tokes
  equations,  \emph{Nonlinear Problems (Proc. Sympos., Madison, Wis.,} Univ. of
  Wisconsin Press, Madison, Wis., 1963, pp.~69--98.

 \bibitem{Si}
J.~Simon,Compact sets in the space {$L\sp p(0,T;B)$},  \emph{Ann. Mat.
  Pura Appl.} (4) \textbf{146} (1987), 65--96.

  
 \bibitem{Sma}
{\v{S}}.~Smagulov, Parabolic approximation of {N}avier-{S}tokes
  equations, \emph{Chisl. Metody Mekh. Sploshn. Sredy} \textbf{10} (1979), no.~1 Gaz.
  Dinamika, 137--149.
  
\bibitem{So95}
C.D.~Sogge, \emph{Lectures on nonlinear wave equations}, Monographs
  in Analysis, II, International Press, Boston, MA, 1995.

\bibitem{Ste93}
E.M.~ Stein, \emph{Harmonic analysis: real-variable methods, orthogonality,
  and oscillatory integrals}, Princeton Mathematical Series, vol.~43, Princeton
  University Press, Princeton, NJ, 1993, With the assistance of Timothy S.
  Murphy, Monographs in Harmonic Analysis, III.

\bibitem{S77}
R. S.~ Strichartz,Restrictions of {F}ourier transforms to quadratic
  surfaces and decay of solutions of wave equations,  \emph{Duke Math. J}. \textbf{44}
  (1977), no.~3, 705--714.
  

\bibitem{Tem69a}
R.~T{\'e}mam, Sur l'approximation de la solution des \'equations de
  {N}avier-{S}tokes par la m\'ethode des pas fractionnaires. {I},  \emph{Arch.
  Rational Mech. Anal}. \textbf{32} (1969), 135--153.

\bibitem{Tem69b}
R.~T{\'e}mam, Sur l'approximation de la solution des \'equations de
  {N}avier-{S}tokes par la m\'ethode des pas fractionnaires. {II}, \emph{Arch.
  Rational Mech. Anal}. \textbf{33} (1969), 377--385.

\bibitem{Tem01}
R.~Temam, \emph{Navier-{S}tokes equations. Theory and numerical analysis}, AMS Chelsea Publishing,
  Providence, RI, 2001, Reprint of the 1984
  edition.
  
\bibitem{YKS}
N.~N. Yanenko, B.~G. Kuznetsov, and Sh. Smagulov, On the approximation of
  the {N}avier-{S}tokes equations for an incompressible fluid by
  evolutionary-type equations, \emph{Numerical methods in fluid dynamics}, ``Mir'',
  Moscow, 1984, pp.~290--314.
  

\end{thebibliography}
\end{document}